\documentclass[12pt]{article}
\usepackage{amsmath,amssymb,amsthm}
\usepackage{graphicx,psfrag,epsf}
\usepackage{dsfont}
\usepackage{enumerate}
\usepackage{natbib}
\usepackage{caption}
\usepackage{subcaption}
\usepackage{floatrow}
\usepackage[hidelinks]{hyperref}

\newcommand{\blind}{0}
\date{\small This preprint has been published in \emph{PeerJ Computer Science}:\\
Durrande N, Hensman J, Rattray M, Lawrence ND. (2016) Detecting periodicities with Gaussian processes. PeerJ Computer Science 2:e50 https://doi.org/10.7717/peerj-cs.50}

\renewcommand{\exp}[1]{\mathrm{exp} \left( #1 \right)}
\newcommand\PS[1]{{\left \langle #1 \right \rangle}_{\! \mathcal{H}}}
\newcommand\N[1]{{\left| \left| #1 \right| \right|}_{\! \mathcal{H}}}
\newcommand\PSi[2]{{ \left \langle #1 \right \rangle}_{\! #2}}

\newcommand\dx{\mathrm{d}}
\def\independenT#1#2{\mathrel{\rlap{$#1#2\mathsurround0pt$}\mkern2mu{#1#2}}}
\newcommand\independent{\protect\mathpalette{\protect\independenT}{\perp}}
\makeatletter
\newcommand{\shorteq}{\ \settowidth{\@tempdima}{a}\resizebox{\@tempdima}{\height}{=}\ }
\makeatother

\DeclareMathOperator*{\Var}{Var}
\DeclareMathOperator*{\E}{E}
\DeclareMathOperator*{\Cov}{Cov}

\begin{document}

\def\spacingset#1{\renewcommand{\baselinestretch}%
{#1}\small\normalsize} \spacingset{1}

%%%%%%%%%%%%%%%%%%%%%%%%%%%%%%%%%%%%%%%%%%%%%%%%%%%%%%%%%%%%%%%%%%%%%%%%%%%%%%

\if0\blind
{
  \title{\bf Gaussian process models for periodicity detection}
  \author{N. Durrande$^1$\thanks{Corresponding author: n.durrande@sheffield.ac.uk}
    , J. Hensman$^1$,  M. Rattray$^2$,
    N. D. Lawrence$^1$ \hspace{.2cm}\\
    $^1$ Department of Computer Science and Sheffield Institute for\\
    Translational Neuroscience, University of Sheffield, UK \hspace{.2cm} \\
    $^2$ Faculty of Life Sciences, University of Manchester, UK\\
}
  \maketitle
} \fi

\bigskip

\hrule
\vspace{0.5cm}
\noindent
\textbf{Summary.}
We consider the problem of detecting the periodic part of a function
given the observations of some input/output tuples $(x_i,y_i)$, $1
\leq i \leq n$. As they are known for being powerful tools for dealing
with such data, our approach is based on Gaussian process regression models which are closely related to reproducing kernel Hilbert
spaces (RKHS). The latter offer a powerful framework for decomposing covariance functions as the sum of periodic and aperiodic kernels. This
decomposition allows for the creation of
sub-models which capture the periodic nature of the signal and its
complement. To quantify the periodicity of the signal, we derive a periodicity
ratio which reflects the uncertainty in the fitted sub-models. Although the method can be
applied to many kernels, we give a special emphasis to the Mat\'ern
family, from the expression of the RKHS inner
product to the implementation of the associated periodic kernels in a
Gaussian process toolkit. The efficiency of the proposed method is
finally illustrated on a biological case study where we detect
periodically expressed genes.\\

\vspace{0.5cm}
\hrule

\spacingset{1.2}
\section{Introduction}
\paragraph{}
The periodic behaviour of natural phenomena arises at many scales, from
the small wavelength of electromagnetic radiations to the movements of
planets. The mathematical study of natural cycles can be traced back
to the XIX century with Thompson's harmonic analysis for predicting
tides~\citep{thomson1878harmonic} and Schuster's investigations on the
periodicity of sunspots~\citep{schuster1898investigation}. Amongst the
methods that have been considered for detecting and extracting the
periodic trend, one can cite harmonic
analysis~\citep{hartley1949tests}, folding
methods~\citep{Stellingwerf1978,leahy1983searches} which are mostly
used in astrophysics and periodic autoregressive
models~\citep{troutman1979some,vecchia1985maximum}. In this article,
we will focus on the application of harmonic analysis in reproducing
kernel Hilbert spaces (RKHS) and on the consequences for Gaussian Process (GP) modelling.

\paragraph{}
Harmonic analysis is based on the projection of a function on a basis of periodic functions. For example, a natural method for extracting the $2\pi$-periodic trend of a function $f$ is to decompose it in a Fourier series:
\begin{equation}
f(x) \rightarrow f_p(x) = a_1 \sin(x) + a_2 \cos(x) + a_3 \sin(2x) + a_4 \cos(2x) + \dots
\end{equation}
where the coefficients $a_i$ are given, up to a
normalising constant, by the $L^2$ inner product between $f$ and the
elements of the basis. However, the phenomenon under study is often
observed in a limited number of points, which means that the value of
$f(x)$ is not known for all $x$ but only for a small set of inputs
$\{x_1, \dots, x_n \}$ called the observation points. With this limited knowledge of $f$, it is
 not possible to compute the integrals of the $L^2$ inner
 product so the coefficients $a_i$ cannot be obtained directly.

\paragraph{}
A popular approach to overcome the fact that $f$ is partially known is
to build a mathematical model $m$ to approximate it. A good model $m$
has to take into account as much information as possible about
$f$. Typically, it interpolates $f$ for the set of observation points
$m(x_i) = f(x_i)$ and its differentiability corresponds to the
assumptions one can have about the regularity of $f$. The main body of
literature tackling the issue of interpolating spatial data is
scattered over three fields:
(geo-)statistics~\citep{matheron1963principles,stein1999interpolation},
functional analysis~\citep{aronszajn1950,berlinet2004reproducing} and
machine learning~\citep{Rasmussen2006}. In the first case, the solution of the interpolation corresponds to the conditional expectation of a Gaussian process $Z$ and in the second, it is the interpolator with minimal norm in a particular Hilbert space $\mathcal{H}$. As many authors pointed out (see for example~\cite{berlinet2004reproducing,Scheuerer2010}), the two approaches are closely related. Both $Z$ and $\mathcal{
H}$ are based on a common object which is a positive definite function
of two variables $k(.,.)$. In statistics, $k$ corresponds to the
covariance of $Z$ and for the functional counterpart, $k$ is the
reproducing kernel of $\mathcal{H}$. From the interpolation or
regularization point of view, the two approaches are equivalent since
they lead to the same model $m$~\citep{wahba1990spline}. Although we
will focus hereafter on the RKHS framework to design periodic kernels, we will also take
advantage of the powerful probabilistic interpretation offered by
Gaussian processes. 

\paragraph{}
A naive approach for extracting the periodic part of $f$ given some
observations would be to approximate it with a mathematical model $m$
and to compute the Fourier coefficients of $m$. However, this method
is not fully satisfactory since each step involves an orthogonal
projection for a different norm. In other words, the construction of
$m$ and the computation of the coefficients are optimal, but not for
the same criterion. As a result the periodic part obtained with this method cannot naturally be seen as a ``best predictor''. To overcome this issue, we propose in this article to build the Fourier series using the RKHS inner product instead of the $L^2$ one. To do so, we extract the sub-RKHS $\mathcal{H}_p$ of periodic functions in $\mathcal{H}$ and model the periodic part of $f$ by its orthogonal projection onto $\mathcal{H}_p$. The prediction then inherits from the probabilistic framework associated with RKHS and the percentage of periodicity of $f$ can elegantly be estimated.

\paragraph{}
The last part of this introduction, gives an overview of the RKHS framework and
emphasises the properties of the Mat\'ern family of kernels. In
section 2, we focus on the construction of periodic kernels. Section 3 details the decomposition of GP models into periodic and aperiodic sub-model. These results allow us to introduce, in Section 4, a new criterion for measuring the periodicity of a signal. Finally, the last section illustrates the proposed approach on a biological case study where we detect, amongst the entire genome, the genes presenting a cyclic expression. This issue of
detecting periodically expressed genes is the application that initially motivated
the present work.

\paragraph{}
The examples and the results presented in this article have been
generated with the version 0.2 of the python Gaussian process toolbox
\textit{GPy}. This toolbox, in which we have implemented the periodic kernels discussed
here, can be downloaded at 
\url{http://github.com/SheffieldML/GPy}.

\subsection{Approximation in Reproducing Kernel Hilbert Spaces}
\paragraph{}
The aim of this section is to introduce the notion of RKHS and to
derive the expression of the best predictor. We will also briefly show
how to construct a RKHS from any positive definite function. For a
more details, we refer the reader
to~\citet[chap. 1]{berlinet2004reproducing} and \cite{aronszajn1950}.

\paragraph{}
Let $\mathcal{H}$ be a Hilbert space of real valued functions defined over $D \subset \mathds{R}$. $\mathcal{H}$ is said to be a RKHS if and only if there exist a function $k(.,.) : D \times D \rightarrow \mathds{R}$ such that for all $x \in D$
\begin{itemize}
 \item[$(i)$] $k(x,.) \in \mathcal{H}$
 \item[$(ii)$] $\forall f \in \mathcal{H},\ f(x)=\PS{f,k(x,.)}$.
\end{itemize}
The function $k$ satisfying these properties is unique and it is called the \textit{reproducing kernel} of $\mathcal{H}$. 

\paragraph{}
Recalling that a function $k$ is said to be positive semi-definite if $\forall m \in \mathds{N},\ \forall \mathbf{a} \in \mathds{R}^m,\  \forall \mathbf{x} \in D^m$
\begin{equation}
\sum_{i=1}^m \sum_{j=1}^m a_i a_j k(x_i,x_j) \geq 0,
 \label{eq:pdkernel}
\end{equation}
it can be shown that a reproducing kernel is necessarily a symmetric
positive semi-definite (spd) function. Reciprocally, the
Moore-Aronszajn theorem states that for all spd-function $k$ on $D
\times D$, there exist only one RKHS of functions on $D$ with $k$ as
reproducing kernel. A common approach is then to define a RKHS by specifying its reproducing kernel. As the covariance of a random process is also a spd-function, we will use interchangeably the words kernel, covariance function and reproducing kernel.

\paragraph{}
To get an insight on the elements of the RKHS associated with a spd-function $k$, we first consider the space $H$ generated by finite combinations of $k(x_i,.)$:
\begin{equation}
H = \left\{ \sum_{i=1}^m a_i k(x_i,.),\ a_i \in \mathds{R},\ x_i \in D,\ m \in \mathds{N} \right\}.
\label{eq:preH}
\end{equation}
Obviously, we have $k(x,.) \in H$ for all $x \in D$ so $(i)$ is
satisfied. Using the property that $k$ is a spd-function, it is straightforward to show that 
\begin{equation}
\PSi{\sum_{i=1}^m a_i k(x_i,.),\sum_{j=1}^{m'} b_j k(x_j,.)}{H} = 
\sum_{i=1}^m \sum_{j=1}^{m'} a_i b_j k(x_i,x_j).
\label{eq:IPH}
\end{equation}
defines a valid inner product on $H$. One particular asset of this inner product is that $k(x,.)$ satisfies $(ii)$. Indeed, for all $ f = \sum a_i k(x_i,.) \in \mathcal{H}$ and $x \in D$ we have
\begin{equation}
\PSi{f,k(x,.)}{H} = 
\sum_{i=1}^m a_i \PSi{k(x_i,.),k(x,.)}{H} = 
\sum_{i=1}^m a_i k(x_i,x) = 
f(x).
\label{eq:repH}
\end{equation} 
Although the properties $(i)$ and $(ii)$ are fulfilled, $H$ is not a
necessarily a RKHS since it may not a Hilbert space (it is not always complete). Let $\mathcal{H}$ be the closure of $H$ and $\PS{.,.}$ the continuous extension of $\PSi{.,.}{H}$ onto $\mathcal{H}$. Then $\mathcal{H}$ is a Hilbert space and it can be shown that $(i)$ and $(ii)$ are still satisfied: $\mathcal{H}$ is the only RKHS with kernel $k$. 

\paragraph{}
We will now focus on how to take advantage of the RKHS framework to
approximate a function $f$ that is observed in a limited number of
points. Let $X=\{x_1,\dots,x_n\} \in D^n$ be a set of points where the
value $y_i=f(x_i)$ is known and $\mathbf{y}$ be the vector of $y_i$. For a given RKHS $\mathcal{H}$, the best interpolator $m$ is defined as the interpolator with minimal norm:
\begin{equation}
 m = \operatornamewithlimits{argmin}_{h \in \mathcal{H}}  \big( \N{h} \
 \big| \ h(x_i)= y_i ,\ i \in 1,\dots,n \big).
\end{equation}
It can be shown that $m$ corresponds to the orthogonal projection of $f$ onto the space spanned by the $k(x_i,.)$:
\begin{equation}
 \mathcal{H}_X = \operatorname{span} \big( k(x_i,.),\ x_i \in X \big).
\end{equation}
Let $\mathbf{k}(.)$ be the $n \times 1$ vector of functions with general term $(\mathbf{k}(.))_i = k(x_i,.)$. This vector corresponds to a basis of $\mathcal{H}_X$. The Gram matrix $\mathbf{K}$ associated to this basis has general term $\mathbf{K}_{ij} = \PS{k(x_i,.),k(x_j,.)} = k(x_i,x_j)$. When $\mathbf{K}$ is invertible, it is straightforward to show that 
\begin{equation}
 k_X(x,y) = \mathbf{k}^T(x) \mathbf{K}^{-1} \mathbf{k}(y) 
\label{eq:KX}
\end{equation}
satisfies $(i)$ and $(ii)$. Since $\mathcal{H}_X$ is a finite
dimensional space it is necessarily complete so $\mathcal{H}_X$ is a
RKHS with reproducing kernel $k_X$. The orthogonal projection of $f$ onto $\mathcal{H}_X$ is then:
\begin{equation}
 m(x) = \PS{k_X(x,.),f} = \mathbf{k}^T(x) \mathbf{K}^{-1} \PS{\mathbf{k}(.),f} = \mathbf{k}^T(x) \mathbf{K}^{-1} \mathbf{y}.
\label{eq:BP}
\end{equation}
In the geostatistical community, $m$ is referred to as the Kriging mean. In the probabilistic framework, this expression corresponds to the
conditional expectation of a centred Gaussian process $Z$ with
covariance $k$ knowing the observations. Furthermore, GP provide naturally some prediction variance for the model:
\begin{equation}
\begin{split}
m(x) & = \E [Z(x) | Z(x_i) \shorteq y_i] = \mathbf{k}^T(x) \mathbf{K}^{-1} \mathbf{y} \\
v(x) & = \Var [Z(x) | Z(x_i) \shorteq y_i] = k(x,x) - \mathbf{k}^T(x) \mathbf{K}^{-1} \mathbf{k}(x)
\end{split}
\label{eq:VP}
\end{equation}

\paragraph{}
One particular asset of Eqs.~\ref{eq:BP}-\ref{eq:VP} is that the expressions of $m,\ v$
only depends on $k$. As a result it is not necessary to derive the
expression of the inner product generated by $k$ to obtain the best
predictor and any spd-function can be used directly to build models.
However, a direct proof of the positive definiteness of a function is
often intractable and a widespread approach is to use well known spd-functions such as the squared-exponential (i.e.\ Gaussian and radial
basis function) or the spline kernel. The next section recalls some
results about another interesting class of spd-functions: the Mat\'ern family. 

\subsection{The Mat\'ern class of kernels}
\paragraph{}
Mat\'ern kernels $k$ are stationary spd-functions, which means that they only depend on the distance between the points they are evaluated at: $k(x,y)=\tilde{k}(|x-y|)$. They are often introduced by the spectral density of $\tilde{k}$~\citep{stein1999interpolation}:
\begin{equation}
S(\omega) =  \left( \frac{ \Gamma (\nu) \theta^{2\nu}}{2 \sigma^2 \sqrt{\pi} \Gamma (\nu + 1/2) (2 \nu)^\nu} \left( \frac{2 \nu}{\theta ^2} + \omega ^2 \right)^{\nu+1/2} \right)^{-1}.
\label{eq:SDmat}
\end{equation}
Three parameters can be identified in this equation: $\nu$ which tunes the differentiability of $\tilde{k}$, $\theta$ which corresponds to a lengthscale parameter and  $\sigma^2$ that is homogeneous to a variance. Note that all these parameters are positive reals. 
\paragraph{}
The actual expressions of the Mat\'ern kernels are simple when the parameter $\nu$ is half-integer. For $\nu=1/2,\;3/2,\;5/2$ we have
\begin{equation}
\begin{split}
k_{1/2}(x,y) & = \sigma^2  \exp{ - \frac{|x-y|}{\theta}} \\
k_{3/2}(x,y) & = \sigma^2 \left( 1+ \frac{\sqrt{3}|x-y| }{\theta} \right) \exp{ - \frac{\sqrt{3}|x-y| }{\theta}} \\
k_{5/2}(x,y) & = \sigma^2 \left( 1+ \frac{\sqrt{5}|x-y| }{\theta} + \frac{5|x-y|^2 }{3 \theta^2} \right) \exp{ - \frac{\sqrt{5}|x-y| }{\theta}}.
\end{split}
\label{eq:MATnu135}
\end{equation}
It can be seen that the parameters $\theta$ and $\sigma^2$ respectively
correspond to a rescaling of the abscissa and ordinate axis. For
$\nu=1/2$ one can recognise the expression of the exponential kernel
(i.e. the covariance of the Ornstein-Uhlenbeck process) and the limit case $\nu \rightarrow \infty $ corresponds to the squared exponential covariance function~\citep{Rasmussen2006}. 

\paragraph{}
One considerable asset of the Mat\'ern class of kernels is to have strong connections with various fields. For example, a Gaussian process $Z$ with Mat\'ern covariance is an autoregressive process. As detailed in appendix~\ref{sec:NORM}, this connection allows to use previous results from the literature to derive the expression of the inner products of the associated RKHS:\\

\noindent
Mat\'ern $1/2$ (exponential kernel) \\
\begin{equation}
\begin{split}
\PSi{g,h}{\mathcal{H}_{1/2}} & = \frac{\theta}{2\sigma^2} \int_a^b \left( \frac{1}{\theta} g + g' \right) \left( \frac{1}{\theta}h + h' \right) \dx t + \frac{1}{\sigma^2} g(a)h(a)
\end{split}
\label{eq:PSk12}
\end{equation}

\noindent
Mat\'ern $3/2$
\begin{equation}
\begin{split}
\PSi{g,h}{\mathcal{H}_{3/2}} & = \frac{\theta^3}{12 \sqrt{3} \sigma^2} \int_a^b \left( \frac{3}{\theta^2} g + 2 \frac{\sqrt{3}}{\theta}g' + g'' \right) \left( \frac{3}{\theta^2} h + 2 \frac{\sqrt{3}}{\theta}h' + h'' \right) \dx t  \\
 & \qquad + \frac{1}{\sigma^2} g(a)h(a) + \frac{\theta^2}{3 \sigma^2} g'(a)h'(a) 
\end{split}
\label{eq:PSk32}
\end{equation}

\noindent
Mat\'ern $5/2$
\begin{equation}
\begin{split}
L_t(g) & =  \sqrt{\frac{3 \theta^5}{400 \sqrt{5} \sigma^2}} \left( \frac{5 \sqrt{5}}{\theta^3} g(t) + \frac{15}{\theta^2}g'(t) + \frac{3 \sqrt{5}}{\theta} g''(t) + g'''(t) \right) \\
\PSi{g,h}{\mathcal{H}_{5/2}} & = \int_a^b L_t (g)L_t (h)  \dx t + \frac{9}{8\sigma^2} g(a)h(a) +  \frac{9 \theta^4}{200 \sigma^2} g(a)''h''(a)  \\
 & \qquad  + \frac{3 \theta^2}{5 \sigma^2} \left(g'(a)h'(a) +\frac18 g''(a)h(a) + \frac18 g(a)h''(a) \right)  
\end{split}
\label{eq:PSk52}
\end{equation}
Although these expressions are direct consequences of \cite{Doob1953}
and \cite{hajek1962linear} they cannot be found in the literature to the best of our knowledge.
 
\paragraph{}
 Another field that is closely related to Mat\'ern kernels is Sobolev
 spaces. As stated in \citet[Theorem 9.1]{porcu2012some} and
 \cite{wendland2005scattered}, the RKHS generated by $k$ coincides
 with the Sobolev space $W^{\nu + 1/2}_2$. This will be particularly
 useful in the next section to show that sine and cosine functions belong
 to the RKHS.

\paragraph{}
\cite{Scheuerer2010} point out that Sobolev spaces are intuitively
more accessible than RKHS (it is often straightforward to tell if a
function belongs or not to $W^n_2$) but RKHS offer a good framework
for deriving an approximation of $f$ based on the observations
$f(x_i)$. As a consequence, Mat\'ern RKHS are very interesting for modelling since they benefit from both assets: the Sobolev structure of $\mathcal{H}$ allows to understand the assumptions on $f$ (for example, $\nu$ is directly linked to differentiability of $f$) and the RKHS properties give a compact expression for the optimal predictor. 

%%%%%%%%%%%%%%%%%%%%%%%%%%%%%%%%%%%%%%%%%%%%%%%%%%%%%%%%%%%%%%%%%%%%%%%%%%%%%%%%
\section{Kernels of periodic subspaces}
%%%%%%%%%%%%%%%%%%%%%%%%%%%%%%%%%%%%%%%%%%%%%%%%%%%%%%%%%%%%%%%%%%%%%%%%%%%%%%%%

%%%%%%%%%%%%%%%%%%%%%%%%%%%%%%%%%%%%%%%%
\subsection{Fourier basis in RKHS}
%%%%%%%%%%%%%%%%%%%%%%%%%%%%%%%%%%%%%%%%

\paragraph{}
In this section, we will see how to extract the subspace of $2
\pi$-periodic functions in a RKHS $\mathcal{H}$. We will assume here
that $\mathcal{H}$ has a Mat\'ern kernel where $\nu$ is
half-integer. However, the method presented here can be applied to any
RKHS as long as the Gram matrix associated to a periodic basis can be
computed. For a detailed list of RKHS inner products we refer the reader to~\cite[Chap. 7]{berlinet2004reproducing}.

\paragraph{}
One popular basis for a space of periodic functions is the Fourier basis $(\sin(x)$, $\cos(x)$, $\sin(2x)$, $\cos(2x)$, $\dots)$. Hereafter, we consider a truncated version of this basis, ignoring the frequencies higher than $q$
\begin{equation}
 \mathbf{F}(x) = (\sin(x), \cos(x), \dots, \sin(qx), \cos(qx))^T,
\end{equation}
and we denote by $\mathcal{H}_p$ the space spanned by this basis. The
fact that $\mathcal{H}$ coincides with $W^{\nu + 1/2}_2$ ensures that the elements of $\mathcal{H}$ are the functions such that
\begin{itemize}
 \item the $i^{\mathrm{th}}$ derivatives $(0 \leq i \leq \nu-1/2)$ are absolutely continuous and square integrable,
 \item the $(\nu + 1/2)^{\mathrm{th}}$ derivative is defined almost
   everywhere and is square integrable.
\end{itemize}
As a consequence, we have $\mathcal{H}_p \subset \mathcal{H}$ since
all the functions of the basis are infinitely differentiable. Let
$\mathbf{G}$ be the Gram matrix of $\mathbf{F}$ in $\mathcal{H}$:
$\mathbf{G}_{i,j} = \PS{\mathbf{F}_i,\mathbf{F}_j}$. Similarly to Eq.~\ref{eq:KX}, it is straightforward to show that
\begin{equation}
 k_{p}(x,y) = \mathbf{F}^T(x) \mathbf{G}^{-1}  \mathbf{F}(y)
\end{equation}
is the reproducing kernel of $\mathcal{H}_p$. Hereafter, we will refer to $k_p$ as the \emph{periodic kernel}.

\paragraph{}
The matrix $\mathbf{G}$ can be computed from the expression of the inner product given in Eqs.~\ref{eq:PSk12}-\ref{eq:PSk52}. In contrast to the Gram matrix of the Fourier basis in $L^2$, $\mathbf{G}$ is not a diagonal matrix if the length of $D$ is a multiple of the period.
One essential property for the practical use of periodic kernels is that the computation of the elements of $\mathbf{G}$ can be performed analytically. Indeed, all the elements of the basis can be written in the form $\cos(\omega x + \varphi)$. Using the notation $L_x$ for the linear operators in the inner product integrals (see Eq.~\ref{eq:PSk52}) we obtain:
\begin{equation}
L_x(\cos(\omega x + \varphi)) = \sum_i \alpha_i \cos(\omega x +
\varphi)^{(i)} = \sum_i \alpha_i \omega^i \cos \left( \omega x +
  \varphi + \frac{i\pi}{2} \right).
\end{equation}
The latter can be factorised in a single cosine $\rho \cos(\omega x + \phi)$ with \\
\begin{equation}
  \rho = \sqrt{  r_c^2 + r_s^2 } \text{, \quad}  \phi = 
\left\{
  \begin{array}{l l}
    \operatorname{arcsin} \left( r_s/\rho \right ) & \quad \text{if }
  r_c \geq 0\\
    \operatorname{arcsin} \left( r_s /\rho \right ) + \pi & \quad\text{if }
  r_c < 0
  \end{array} \right.
\end{equation}
where $ \displaystyle r_c = \sum_i \alpha_i \omega^i \cos \left( \varphi +
  \frac{i\pi}{2} \right)$ and $ \displaystyle r_s = \sum_i \alpha_i \omega^i \sin
\left( \varphi + \frac{i\pi}{2} \right) $.

Eventually, the computation of the inner product boils down to the
integration of a product of two cosines, which can be solved by
linearisation.

\subsection{Tuning the period}
\label{sec:estper}
\paragraph{}
We assumed previously a $2 \pi$-periodicity for the signal. However this period can be modified by introducing a
parameter $\lambda$ in the definition of the Fourier basis:
\begin{equation}
\mathbf{F}_\lambda(x) = \left(\sin \bigg(\frac{2 \pi}{\lambda} x \bigg),\  \dots,\
\cos \bigg(\frac{2 \pi}{\lambda} qx \bigg) \right)^T.
\end{equation}
As for the other parameters of the kernel $\sigma^2$ and $\theta$, maximum likelihood
estimation can be used to obtain a value of $\lambda$ well suited to
the data. This estimation of the period will be illustrated in the
case study of the next section.

%%%%%%%%%%%%%%%%%%%%%%%%%%%%%%%%%%%%%%%%%%%%%
\subsection{Application to a benchmark}
\paragraph{}
We will now illustrate on a benchmark of test functions the use periodic kernels for GP modelling. Furthermore, we will compare the resulting models with COSOPT~\citep{straume2004dna} and ARSER~\citep{yang2010analyzing} which are representative of the methods commonly used in biostatistics for detecting periodically expressed genes~\citep{hughes2009harmonics,amaral2012circadian}. 

\paragraph{}
COSOPT assumes the following model for the signal:
\begin{equation}
y(t) = \alpha + \beta t + \gamma \cos(\omega t + \varphi) + \varepsilon,
\end{equation}
where $\varepsilon$ corresponds to some white noise. The algorithm proceeds in two steps to determine the values of $\alpha$, $\beta$, $\gamma$, $\omega$ and $\varphi$. First, a linear regression model is fitted to estimate the value of $\alpha,\ \beta$. The linear trend is then subtracted from the signal and the remaining parameters are fitted by minimizing the mean square error.

\paragraph{}
The underlying model is more sophisticated for ARSER since it accounts for various frequencies: 
\begin{equation}
y(t) = \alpha + \beta t + \sum_i \gamma_i \cos(\omega_i t + \varphi_i) + \varepsilon.
\end{equation}
The number of cosine terms and their frequencies $\omega_i$ are obtained by detecting the peaks of the spectrum via the fast Fourier transform. In practice, this number is typically around 1-5. As previously the linear trend is initially subtracted to the data and an additional smoothing is performed to limit high frequencies due to noise. Although this model is more flexible, it has two drawbacks: the input points are assumed to be regularly spaced and the resulting model is not necessarily periodic. 

\paragraph{}
In addition, we introduce the following Gaussian process model:
\begin{equation}
Y(t) = \alpha + \beta t + Y_p(t) + \varepsilon.
\label{eq:GPpermod}
\end{equation}
where $Y_p$ has periodic kernel $k_p$. Here, $\alpha$ and $\beta$ should be interpreted as random variables with Gaussian distribution $\mathcal{N}(0,1)$. The best predictor associated with this model given by Eq.~\ref{eq:BP}) where $k(x,y) = 1 + xy + \sigma_p^2 k_p(x,y) + \tau^2 \delta (x,y)$. The parameters $\sigma_p^2$, $\theta$, $\lambda$ and $\tau^2$ are obtained by maximising the likelihood of the observations. This is equivalent to minimizing -2 times the log-likelihood:
\begin{equation}
 \mathcal{L} = n \log (2 \pi) + \log | \mathbf{K} | + \mathbf{Y}^T \mathbf{K}^{-1}\mathbf{Y}
\label{eq:loglik}
\end{equation}
which depends on all the parameters of the kernel through the matrix $\mathbf{K}$. The number of frequencies in the Fourier basis is set to $q=20$. The best predictor associated with this model given by Eq.~\ref{eq:BP} where $k(x,y) = 1 + xy + k_p(x,y) + \tau^2 \delta_{x,y}$. Although this information is readily available, we will not use in this benchmark the prediction variance provided by the GP models.

\paragraph{}
COSOPT and ARSER also include additional features for measuring the periodicity of a signal in term of $p$-value or false discovery rate. The probabilistic framework of Gaussian processes could be used to derive such statistics for the proposed model but these developments fall out of the scope of the present article.

\begin{figure}
    \centering
    \includegraphics[width=0.8\textwidth]{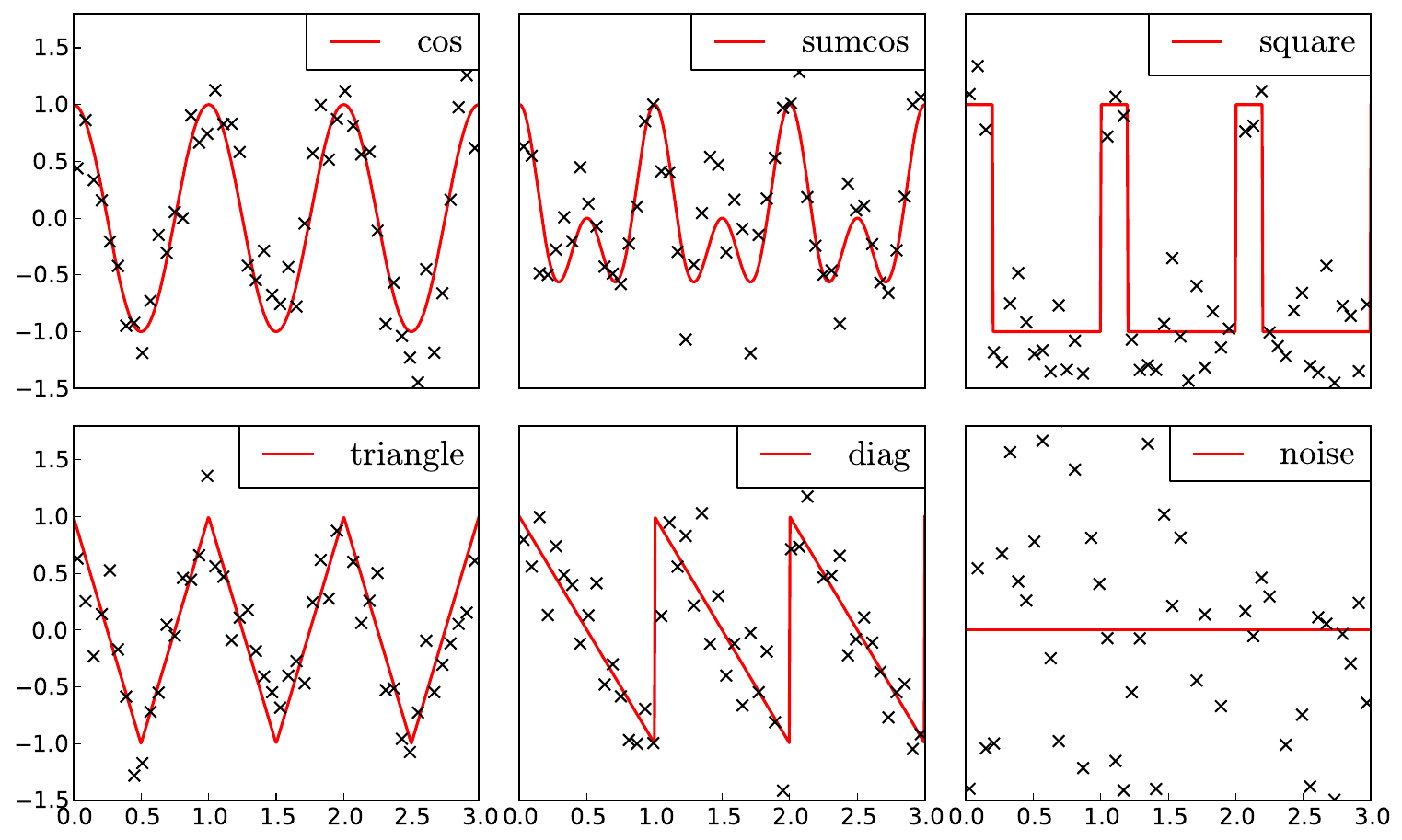}
    \caption{Test functions considered in the benchmark. The crosses indicate the observed values after adding the random noise.}
\label{fig:testfunct}
\end{figure}

\paragraph{}
The prediction of these models are compared on a benchmark of 1-periodic test functions defined over $[0,3]$: $\cos (2 \pi t)$, $\mathrm{sumcos}(t) = 1/2 (\cos (2 \pi t) + \cos(4 \pi t))$, $\mathrm{square}(t)$, $\mathrm{triangle}(t)$,  $\mathrm{diag}(t)$ and $\mathrm{noise}(t)$ which are represented in Figure~\ref{fig:testfunct}. A training set of 50 equally spaced test points is used for learning these functions and a $\mathcal{N}(0,0.1)$ observation noise is added to each observation, except for $\mathrm{noise}$ where the perturbations are $\mathcal{N}(0,1)$. As the test functions do not include any linear trend, the value of $\beta$ is fixed to zero for all models.

\begin{figure}
    \centering
    \includegraphics[width=0.8\textwidth]{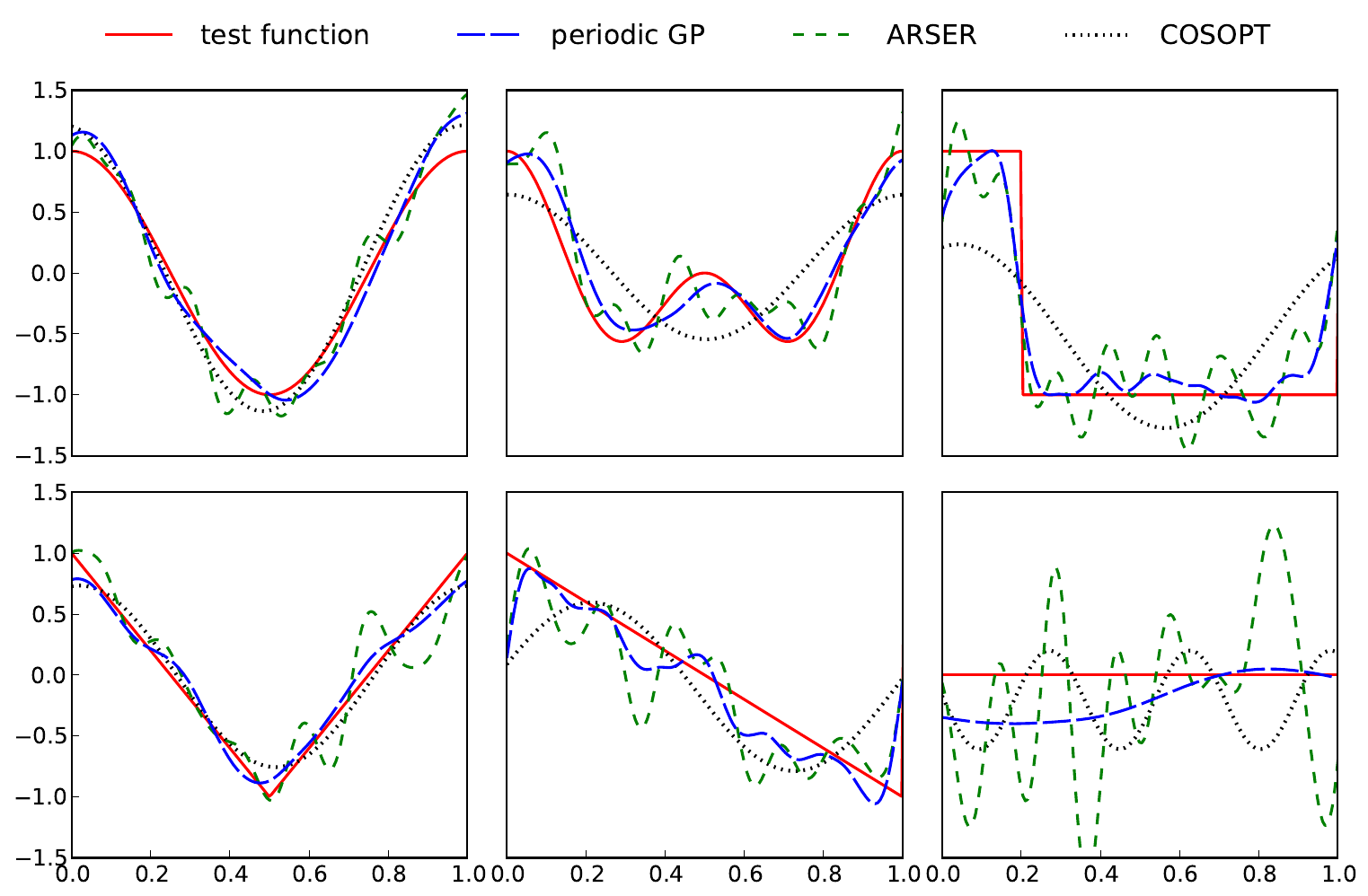}
    \caption{Plots of the test functions with associated fitted models. For an improved visibility, the plotting region is limited to one period. The periodic GP model is based on a periodic Mat\'ern kernel with regularity $\nu = 3/2$.}
\label{fig:compfit}
\end{figure}

\paragraph{}
The models fitted with COSOPT, ARSER and the periodic GP can be compared in Figure~\ref{fig:compfit}. To asses the overall precision, we repeat the fitting procedure 50 times with different values of the observation noise. The root mean square error (RMSE) is computed based on a 500-point test set spanning $[0,3]$. A summary of the obtained result is given in Table~\ref{tab:RMSE}.

\begin{table}
\begin{tabular}{|c|ccccc|}
\hline
test function & COSOPT & ARSER & GP $\nu = 1/2$ & GP $\nu = 3/2$ & GP $\nu = 5/2$ \\ \hline
cos      &  \textbf{\emph{0.09}}  (0.03)  &  0.23  (0.03)  &  0.16  (0.03)  &  \textbf{0.11}  (0.03)  &  \textbf{0.10}  (0.03) \\
sumcos   &  0.36  (0.01)  &  0.24  (0.08)  &  0.19  (0.09)  &  \textbf{0.14}  (0.04)  &  \textbf{\emph{0.13}}  (0.04) \\
square   &  0.60  (0.01)  &  0.37  (0.03)  &  \textbf{\emph{0.31}}  (0.05)  &  \textbf{0.32}  (0.04)  &  \textbf{0.32}  (0.03) \\
triangle   &  \textbf{\emph{0.11}}  (0.02)  &  0.23  (0.03)  &  0.15  (0.03)  &  \textbf{0.12}  (0.03)  &  \textbf{0.12}  (0.03) \\
diag     &  0.36  (0.01)  &  0.33  (0.04)  &  \textbf{\emph{0.26}}  (0.04)  &  \textbf{\emph{0.26}}  (0.03)  &  \textbf{\emph{0.26}}  (0.03) \\
noise    &  \textbf{0.40}  (0.06)  &  0.73  (0.11)  &  \textbf{0.44}  (0.20)  &  \textbf{0.39}  (0.19)  &  \textbf{\emph{0.37}}  (0.21) \\ \hline
 mean &  0.32  &  0.36  &  0.25  &  0.22  &  0.22 \\ \hline
\end{tabular}
\caption{Mean value (and standard deviation) of RMSE for each test function and model. The best fit is indicated in italic. The models within one standard deviation from the best result are indicated in bold.}
\label{tab:RMSE}
\end{table}

\paragraph{}
Many remarks can be formulated based on the observation of Figure~\ref{fig:compfit} and Table~\ref{tab:RMSE}. First, COSOPT gives a good fit for the cosine function, but also for the $\mathrm{triangular}$ test function. This can be explained by the overall cosine shape of the latter. The noise filtering can be judged satisfactory for this model but, as expected, the model fails to approximate  non-sinusoidal patterns such as $\mathrm{square}$ and $\mathrm{sumcos}$. The wider range of frequencies allowed in ARSER makes it capable of approximating these more complicated patterns. The drawback for this model is its sensitivity to noise. Indeed high frequencies oscillations corresponding to noise overfitting can be observed on ARSER models. Although some functions in the test set are typically difficult to approximate with Gaussian process models due the presence of discontinuity, the models based on periodic kernels perform remarkably well on this benchmark. On the one hand, the large number of frequencies considered in the truncated Fourier basis allows a good fit of non-sinusoidal patterns. On the other, the embedding of this basis into a Mat\'ern RKHS naturally imposes a penalty on the high frequencies which results in a good filtering of the noise. 

\paragraph{}
Note that the results obtained for the periodic GP models are not specific to the class of periodic kernel introduced in this article. Usual periodic kernels such as $k(x,y) = \exp{- (\sin(|x-y|)^2}$ (see \cite{Rasmussen2006} for more details) would probably lead to similar results. We will detail in the next section one particular asset of the proposed kernels in term of decomposition of the signal.

%%%%%%%%%%%%%%%%%%%%%%%%%%%%%%%%%%%%%%%%%%%%%%%%%%%%%%%%%%%%%%%%%%
%%%%%%%%%%%%%%%%%%%%%%%%%%%%%%%%%%%%%%%%%%%%%%%%%%%%%%%%%%%%%%%%%%
\section{Decomposition of models}

\subsection{Decomposition of kernels}
\begin{figure}
    \begin{subfigure}[b]{0.33\textwidth}
            \centering
            \includegraphics[width=\textwidth]{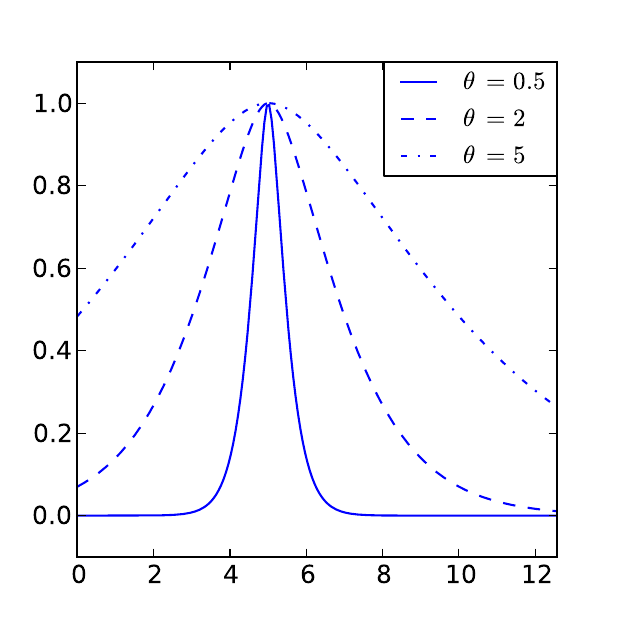}
            \caption{Mat\'ern kernel $k$.}
    \end{subfigure}%
    \begin{subfigure}[b]{0.33\textwidth}
            \centering
            \includegraphics[width=\textwidth]{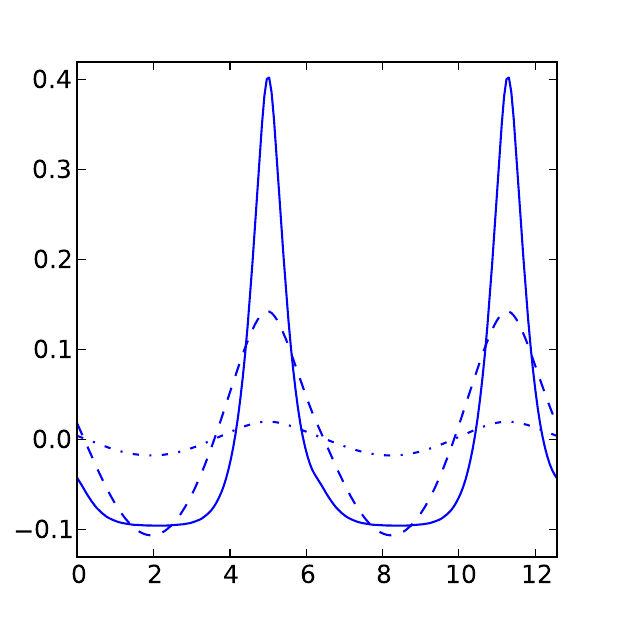}
            \caption{periodic sub-kernel $k_p$.}
    \end{subfigure}
    \begin{subfigure}[b]{0.33\textwidth}
            \centering
            \includegraphics[width=\textwidth]{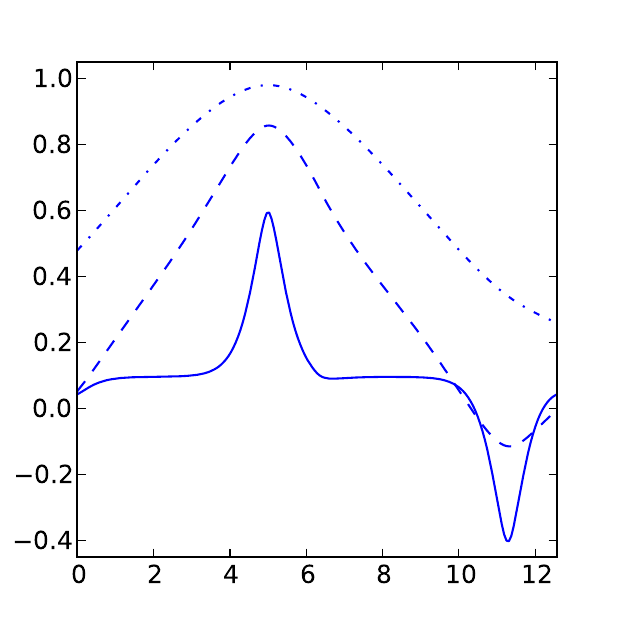}
            \caption{aperiodic sub-kernel $k_a$.}
    \end{subfigure}
    \caption{Examples of decompositions of a Mat\'ern 3/2 kernel
        as a sum of a periodic and aperiodic sub-kernels. The three
        graphs on each plot correspond to a different value of the
        lengthscale parameter $\theta$. For this example the input
        space is $D=[0,4 \pi]$, the cut-off frequency is $q=20$ and one of the
        variables of the kernels is arbitrarily fixed to $5$.}
\label{fig:dec_m32}
\end{figure}

\paragraph{}
The difference of two kernels is generally not a valid covariance function. However the construction of $k_{p}$ ensures that, in this particular case, $k_a = k - k_{p}$ corresponds to a kernel. This is straightforward to see using the RKHS framework since $k_{a}$ is the reproducing kernel of the orthogonal complement of $\mathcal{H}_p$ in $\mathcal{H}$ \citep{berlinet2004reproducing}. As this space is orthogonal to the (truncated) Fourier basis, it will be referred to as the subspace of \emph{aperiodic} functions (hence the subscript $a$). From the probabilistic point of view, this decomposition corresponds to the decomposition of $Z$ as a sum of two independent Gaussian processes, with covariance functions $k_{p}$ and $k_a$. This can be summarized as follow:
\begin{equation}
 \begin{split}
  k = k_{p} + k_a, \qquad \mathcal{H} = \mathcal{H}_p \stackrel{\perp}{+} \mathcal{H}_a, \qquad Z = Z_p \stackrel{\independent}{+} Z_a.
 \end{split} 
\label{eq:decspace}
\end{equation}
An illustration of the decomposition of Mat\'ern 3/2 kernels can be found in Figure~\ref{fig:dec_m32}.

\subsection{Periodic and aperiodic sub-models}
\paragraph{}
The expressions of Eq.~\ref{eq:decspace} allow decomposition of the best predictor as a sum of two sub-models $m_{p}$ and $m_{p}$: 
\begin{equation}
\begin{split}
  m(t) & = \E [Z_p(t) + Z_a(t) | Z(t_i) \shorteq y_i] \\
  & = \E[Z_p(t) | Z(t_i) \shorteq y_i] + \E[Z_a(t) | Z(t_i) \shorteq y_i] \\
  & = \mathbf{k}_{p}(t)^T \mathbf{K}^{-1} \mathbf{y}  + \mathbf{k}_a(t)^T \mathbf{K}^{-1} \mathbf{y}.
\end{split}
\label{eq:decm}
\end{equation}
Similarly, prediction variances are associated to the sub-models
\begin{equation}
\begin{split}
  v_{p}(t) & = \Var[Z_p(t) | Z(t_i) \shorteq y_i] = k_{p}(t,t) - \mathbf{k}_{p}(t)^T \mathbf{K}^{-1} \mathbf{k}_{p}(t) \\
  v_a(t) & = \Var[Z_a(t) | Z(t_i) \shorteq y_i] = k_a(t,t) - \mathbf{k}_a(t)^T \mathbf{K}^{-1} \mathbf{k}_a(t).
\end{split}
\label{eq:decv}
\end{equation}
However, contrarily to Eq.~\ref{eq:decm}, we have $v(t) \neq v_{p}(t) + v_{a}(t)$ since $Y_{p}$ and $Y_{a}$ are not independent knowing the observations. For a
detailed discussion on the decomposition of models based on a sum of
kernels see~\cite{durrande2012additive}.

\paragraph{}
The sub-models can be interpreted as usual GP models with correlated noise. For example, $m_{p}$ is the best predictor based on kernel $k_{p}$ with an observational noise given by $\mathbf{K}_a$. For the RKHS framework, $m_{p}$ and $m_a$ correspond to the solution of a regularization problem and they respectively belong to $\mathcal{H}_p$ and $\mathcal{H}_a$.

\paragraph{}
We now illustrate this model decomposition on the Mauna Loa Observatory dataset \citep{keeling2005atmospheric} which is frequently used in modelling~\citep{Rasmussen2006,wilson2013gaussian}. This dataset contains the monthly average of CO$_2$ concentration in the atmosphere since 1958,  expressed in micromol of CO$_2$ per mol of dry air. Hereafter we will focus on the first six years of the time series, using the initial 48 time points as training data and predicting for the following 24 months. For this dataset we will assume that the one-year period is known.

\begin{figure}
        \begin{subfigure}[b]{0.33\textwidth}
                \centering
                \includegraphics[width=\textwidth]{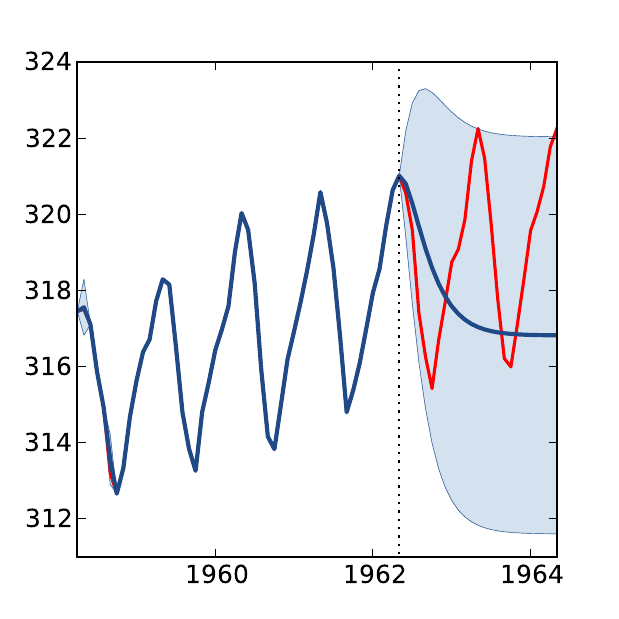}
                \caption{global model $m$.}
        \end{subfigure}%
        \begin{subfigure}[b]{0.33\textwidth}
                \centering
                \includegraphics[width=\textwidth]{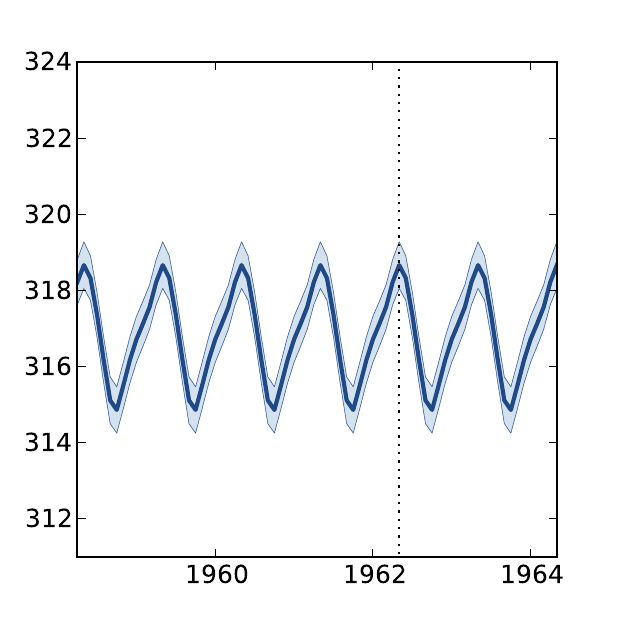}
                \caption{periodic sub-model $m_p$.}
        \end{subfigure}
        \begin{subfigure}[b]{0.33\textwidth}
                \centering
                \includegraphics[width=\textwidth]{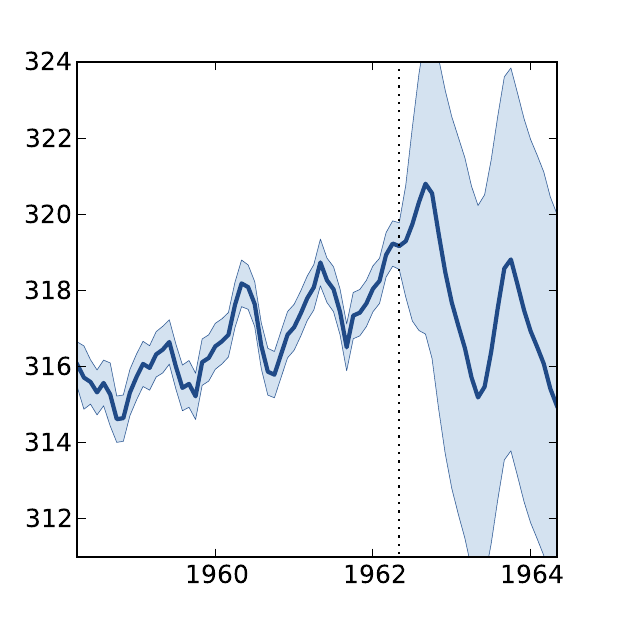}
                \caption{aperiodic sub-model $m_a$.}
        \end{subfigure}

        \caption{Decomposition of a model based on the Mauna Loa Observatory dataset.
          The model is trained on the 48 data-points contained in the left part of the graph. The
          kernel is Mat\'ern 3/2, and the cut-off parameter of the Fourier basis is set to $q=20$. The shaded area corresponds to 95\% confidence intervals and the test function is represented in red. The small increase of the confidence interval width in the left of panel a is due to missing data, which is naturally supported by GP models.}
\label{fig:dec_M52}
\end{figure}

\paragraph{}
 We first consider a GP regression model based on a regular Mat\'ern 3/2 kernel, with maximum likelihood estimation of $\sigma^2$, $\theta$ and $\tau^2$. Figure~\ref{fig:dec_M52} represents the decomposition of the model as detailed in Eqs.~\ref{eq:decm}-\ref{eq:decv}. It can be seen that the periodic sub-model 
successfully extracts the periodic component. Although this model gives very accurate predictions in the training region it drastically fails to forecast the behaviour of the signal after the last observation. We will now see how to improve this result using the sub-kernels.

\subsection{Parametrisation of the kernel}

\paragraph{}
A Mat\'ern kernel $k$ initially depends on three parameters: the
regularity $\nu$, its variance $\sigma^2$ and its lengthscale
$\theta$. However, the decomposition $k = k_{p}
+ k_a$ allows us to set the values of those parameters separately for
each sub-kernel in order to increase the flexibility of the model. The new set of parameters of $k$ is then $(\nu_{p},\ \sigma^2_{p},\ \theta_{p},\ \nu_a,\ \sigma^2_a,\ \theta_a)$, to which $\lambda$ may be added if the period is unknown. 

\paragraph{}
After reparametrisation, $k$ belongs to a larger family of kernels that encapsulates the Mat\'ern one. Furthermore, if $\nu_{p} = \nu_{a}$, and
$\sigma^2_{p},\ \sigma^2_a \neq 0$ the RKHS generated by $k$ and the
one associated with a Mat\'ern kernel with equal regularity correspond to the same space, but endowed with a different norm.

\begin{figure}
        \begin{subfigure}[b]{0.33\textwidth}
                \centering
                \includegraphics[width=\textwidth]{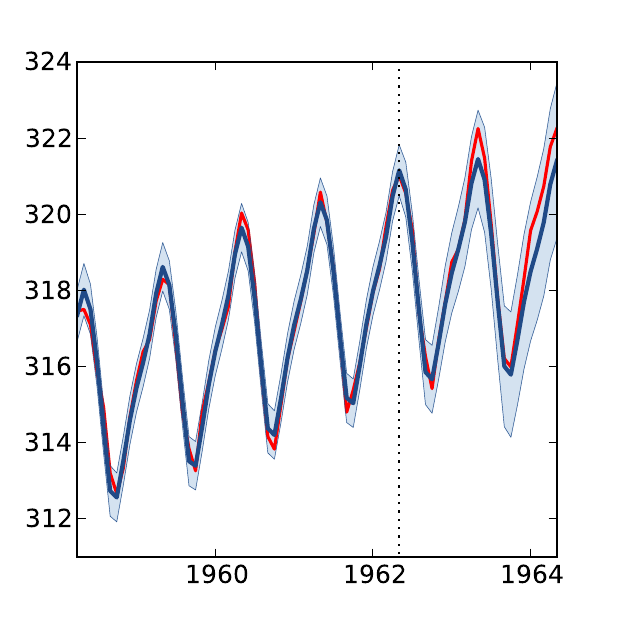}
                \caption{global model $m$.}
        \end{subfigure}%
        \begin{subfigure}[b]{0.33\textwidth}
                \centering
                \includegraphics[width=\textwidth]{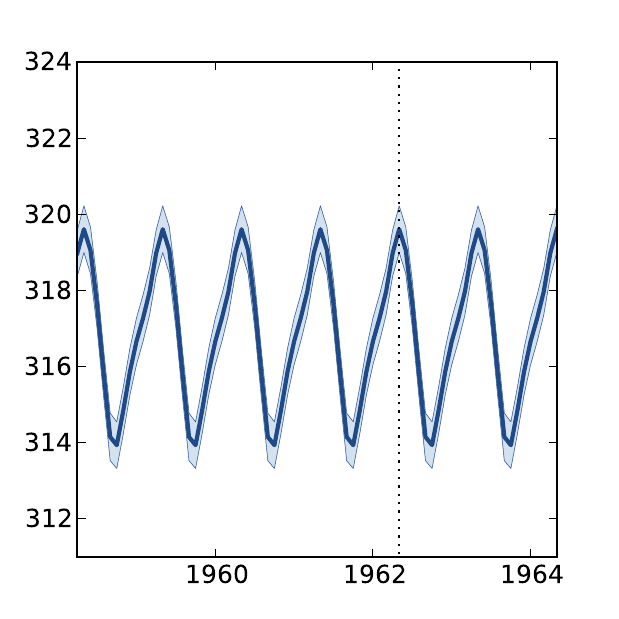}
                \caption{periodic sub-model $m_p$.}
        \end{subfigure}
        \begin{subfigure}[b]{0.33\textwidth}
                \centering
                \includegraphics[width=\textwidth]{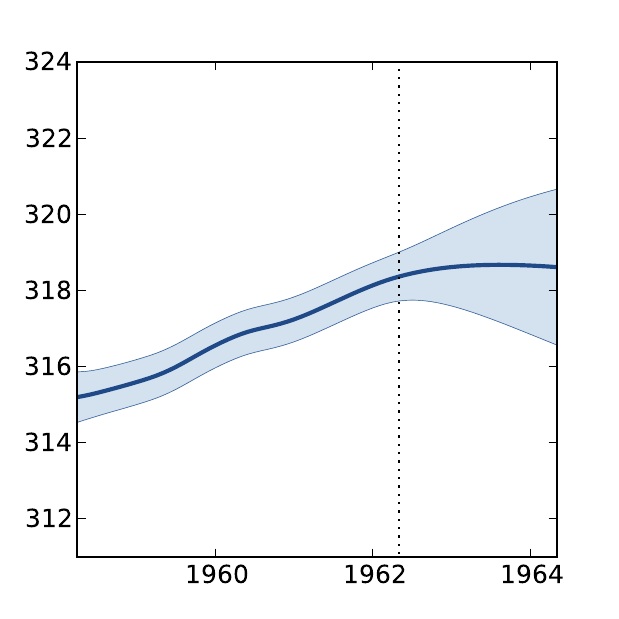}
                \caption{aperiodic sub-model $m_a$.}
        \end{subfigure}
        \caption{Model and sub-models after parametrisation of
          the kernel by $(\sigma^2_{p},\ \theta_{p},\ \sigma^2_a,\
          \theta_a)$. The test points and the other settings are the
          same as in Figure~\ref{fig:dec_M52}.}
\label{fig:dec_M52bis}
\end{figure}

\paragraph{}
The graphs presented in Figure~\ref{fig:dec_M52bis} show the obtained model
after estimating $(\sigma^2_{p},\ \theta_{p},\ \sigma^2_a,\ \theta_a)$, the regularities
$(\nu_{p},\ \nu_a)$ being fixed to $3/2$. In this example, adding two parameters drastically improves the fit of the test function outside the observation region. The global behaviour of the phenomenon is successfully captured by the model which is capable of reproducing both the small scale patterns (oscillations) and the large scale trend. One limitation here is that the regularity parameter of the periodic and aperiodic sub-models is assumed to be the same whereas observation of data suggests a smaller differentiability order for the periodic part.

\section{Measuring the periodicity}
\paragraph{}
The decomposition of the model into a sum of sub-models can be used for
estimating a ratio of periodicity of the signal. In sensitivity
analysis, a common approach for measuring the effect of a set of
variables $(x_1,\dots,x_n)$ on the output of a multivariate function
$f(x_1,\dots,x_n)$ is to introduce a random vector $\mathbf{X}=(X_1,\dots,X_n)$
with values in the input space of $f$ and to define the variance
explained by one subset of variables
$\mathbf{x}_I=(x_{I_1},\dots,x_{I_m})$ as $V_I = \Var \left( \E
  \left( f(\mathbf{X}) | \mathbf{X}_I \right)
\right)$ \citep{oakley2004probabilistic}. Furthermore, the
probabilistic features of the GP model can be taken into account by
computing the indices based on random paths of the conditional GP~\citep{marrel2009calculations}.

\paragraph{}
We now apply these two principles to define a periodicity ratio based
on the sub-models. Let $T$ be a random variable defined over the input
space and $Z_p$, $Z_a$ be the periodic and aperiodic
components of the conditional GP $Z$ knowing it
interpolates the data-points. $Z_p$ and $Z_a$ are normally distributed
with respective mean and variance $(m_{p},\ v_p)$, $(m_{a},\ v_{a})$ and their covariance is
given by $\Cov(Z_p(t),Z_a(t')) = - \mathbf{k}_{p}(t)^T
\mathbf{K}^{-1} \mathbf{k}_{a}(t') $. To quantify the periodicity of the signal we introduce the following periodicity ratio:
\begin{equation}
 R = \frac{\Var_T [Z_p(T)]}{\Var _T [Z_p(T) + Z_a(T)]}.
\label{eq:RatioPer}
\end{equation}
Note that $R$ does not correspond to the percentage of periodicity of the signal in a rigorous way since the dependence between $Z_p$ and $Z_a$ implies $\Var_T [Z(T)] \neq \Var _T [Z_p(T)] + \Var _T [Z_a(T)]$.

\section{Application to gene expression studies}
\paragraph{}
The 24 hour cycle of days can be observed in the oscillations of
biological mechanisms at many scales. This phenomenon, called
circadian rhythm, can for example be seen at a microscopic level on
gene expressions. The cellular mechanism ensuring this periodic
behaviour is called the circadian clock. For \textit{Arabidopsis},
which is a widely used organism in plant biology and genetics, the
study of the circadian clock at a gene level shows an auto-regulatory
system involving several genes~\citep{ding2007complex}. As advocated
in~\cite{edwards2006flowering}, it is believed that the genes involved
in the oscillatory mechanism have a cyclic expression so the detection
of periodically expressed genes is of great interest for completing
current models. As stated in the introduction, this application is the
one that motivated the work presented in this article. 

\paragraph{}
The mechanism allowing genes to interfere in the functioning of the cell
can be summarised as follows: DNA is first duplicated into messenger
RNA, and this RNA is then used for protein synthesis. To quantify the
expression of a specific gene it is thus possible to measure the
concentration of RNA molecules associated with this gene. Microarray
analysis and RNA-sequencing are two examples of methods that take
advantage of this principle. 

\paragraph{}
The dataset we consider here has been initially studied
by~\cite{edwards2006flowering}\footnote{The original dataset is available online at
  \url{http://millar.bio.ed.ac.uk/data.htm}.}. It corresponds to gene
expression for nine day old \textit{arabidopsis} seedlings. After eight days
under a 12h-light/12h-dark cycles, the seedlings are transferred into
constant light. A microarray analysis is performed every four
hours, from 26 to 74 hours after the last dark-light transition, to
monitor the expression of 22810 genes. \cite{edwards2006flowering} use COSOPT~\citep{straume2004dna} for detecting periodicity genes and
identify a subset of 3504 periodically expressed genes, with an estimated a period between 20 and 28 hours. 

\paragraph{}
We now apply to this dataset the method described in the previous
sections. The kernel we consider is a sum of a periodic and aperiodic
Mat\'ern 3/2 kernel plus a delta function to reflect observation noise:
\begin{equation}
 k(t,t') = \sigma_p^2 k_p(t,t') + \sigma^2_a k_a(t,t') + \tau^2 \delta(t,t').
\end{equation}
Although the cycle of the circadian clock is known to be around 24
hours, circadian rhythms often depart from this figure (indeed
circadian is Latin for \emph{around a day}) so we introduce a parameter
$\lambda$ as in Sec.~\ref{sec:estper} to estimate the actual period. The final
parametrisation of $k$ is based on six variables: $(\sigma^2_p,
\theta_p, \sigma^2_a, \theta_a, \tau^2, \lambda)$. For each gene, the
values of these parameters are estimated using maximum likelihood. The
optimization is based on the standard options of the GPy toolkit with
the following boundary limits for the parameters:
$\sigma_p,\ \sigma_a \geq 0$; 
$\theta_p, \ \theta_a \in [10,\ 60] $; $\tau^2 \in
[10^{-5},0.75]$ and $\lambda \in [20,\ 28]$. Furthermore 50 random restarts
are performed for each optimization to limit the effects of local
minimums. 

\paragraph{}
Eventually, the periodicity of each model is assessed with the ratio $R$ given by Eq.~\ref{eq:RatioPer}. As this ratio is a random variable, we approximate the expectation of $R$ with the mean
value of 1000 realisations. To obtain results comparable with the original
paper on this dataset, we label as periodic the set of 3504 genes with
the highest periodicity ratio. The cut-off periodicity ratio associated with
this quantile is 0.77.

\begin{figure}
\begin{floatrow}
\floatbox[]{table}[5cm][][c]{
\begin{tabular}{c|cc}
\# of genes & $\mathcal{P}_{GP}$ & $\overline{\mathcal{P}_{GP}}$ \\ \hline
$\mathcal{P}_{COSOPT}$ & 2127 & 1377 \\
$\overline{\mathcal{P}_{COSOPT}}$ & 1377 & 17929 \\
\end{tabular}
\vspace{1.8cm}
}{
  \caption{Confusion table associated to the predictions by COSOPT and
    the proposed GP approach.}
  \label{tab:compres}
} \qquad
\ffigbox{
\includegraphics[width=6cm]{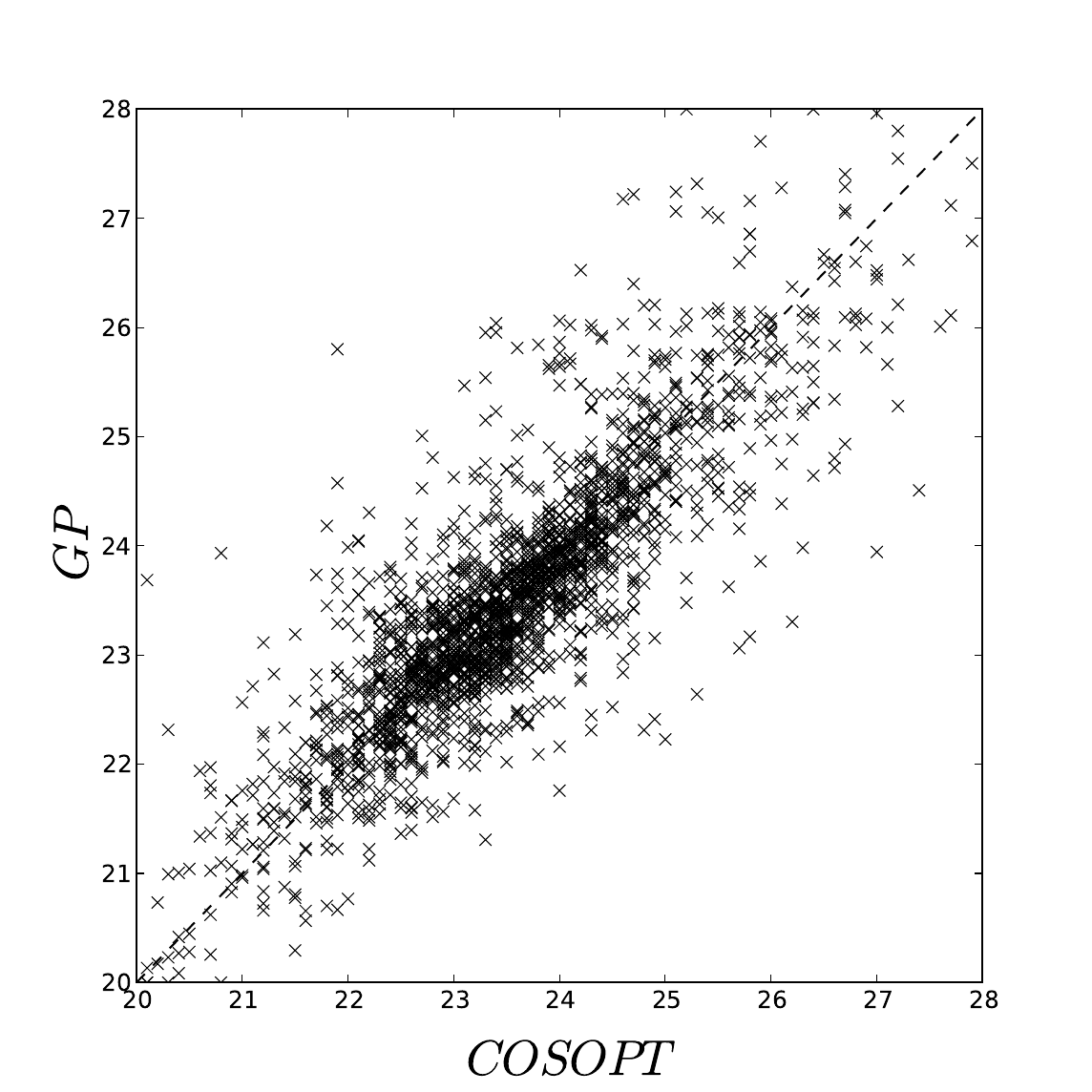}
}{
  \caption{Estimated periods for the genes in $\mathcal{P}_{GP} \cap \mathcal{P}_{COSOPT}$. The coefficient of
  determination of $x \rightarrow x$ (dashed line) is 0.69.}
  \label{fig:estper}
}
\end{floatrow}
\end{figure}

\paragraph{}
Let $\mathcal{P}_{COSOPT}$ and $\mathcal{P}_{GP}$ be the sets
of selected periodic genes respectively by~\cite{edwards2006flowering} and the
method presented here. The overlap between the two sets is summarised
in Table~\ref{tab:compres} where $\overline{\mathcal{S}}$ denotes the
complement of a subset $\mathcal{S}$. Although the results cannot be compare to
any ground truth, the methods seem coherent since 88\% of the genes
share the same label. Furthermore the estimated value of the period
$\lambda$ is consistent for the genes labelled as periodic
by the two methods, as seen in Figure~\ref{fig:estper}.

\begin{figure}
\begin{center}
        \begin{subfigure}[b]{0.45\textwidth}
                \includegraphics[width=\textwidth]{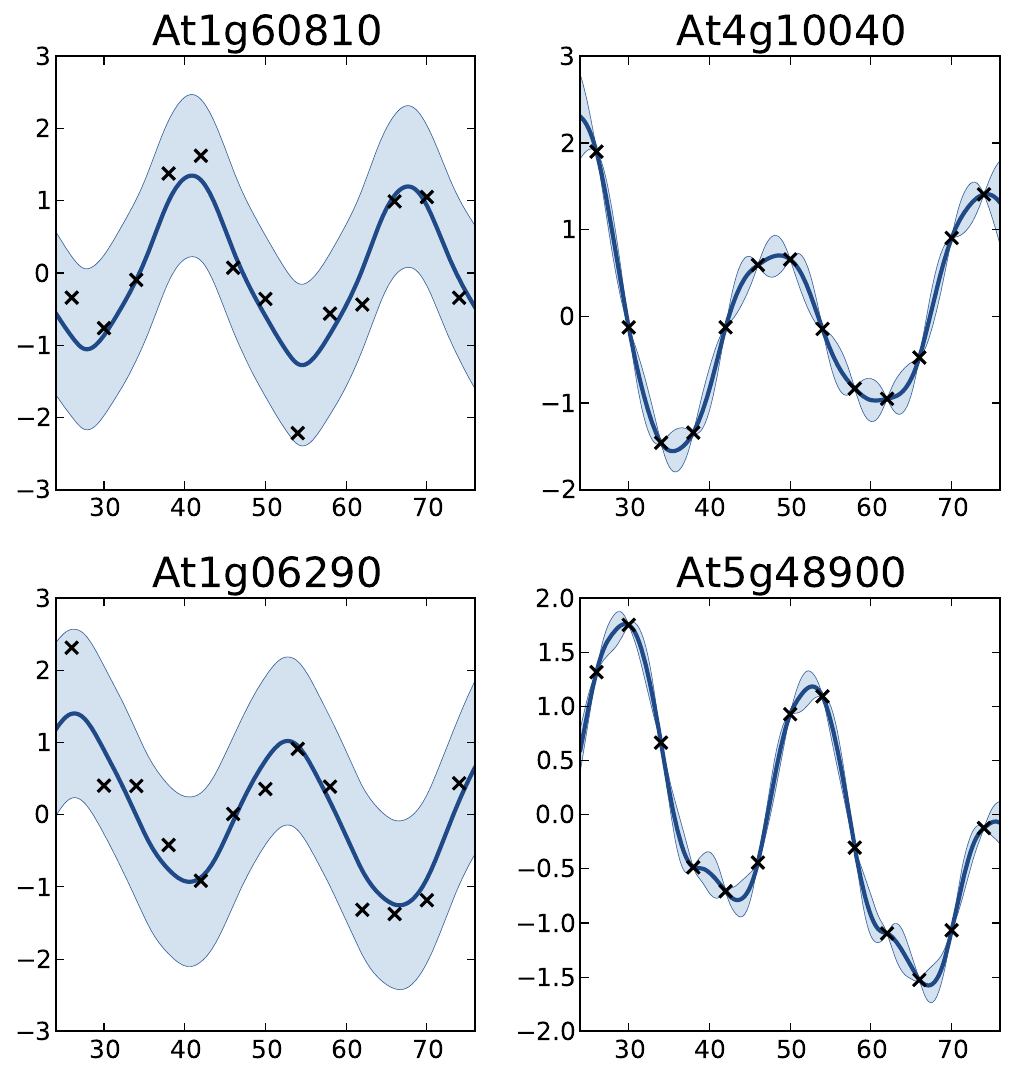}
                \caption{Genes labelled as periodic only by COSOPT.}
        \end{subfigure} \hspace{1cm}
        \begin{subfigure}[b]{0.45\textwidth}
                \centering
                \includegraphics[width=\textwidth]{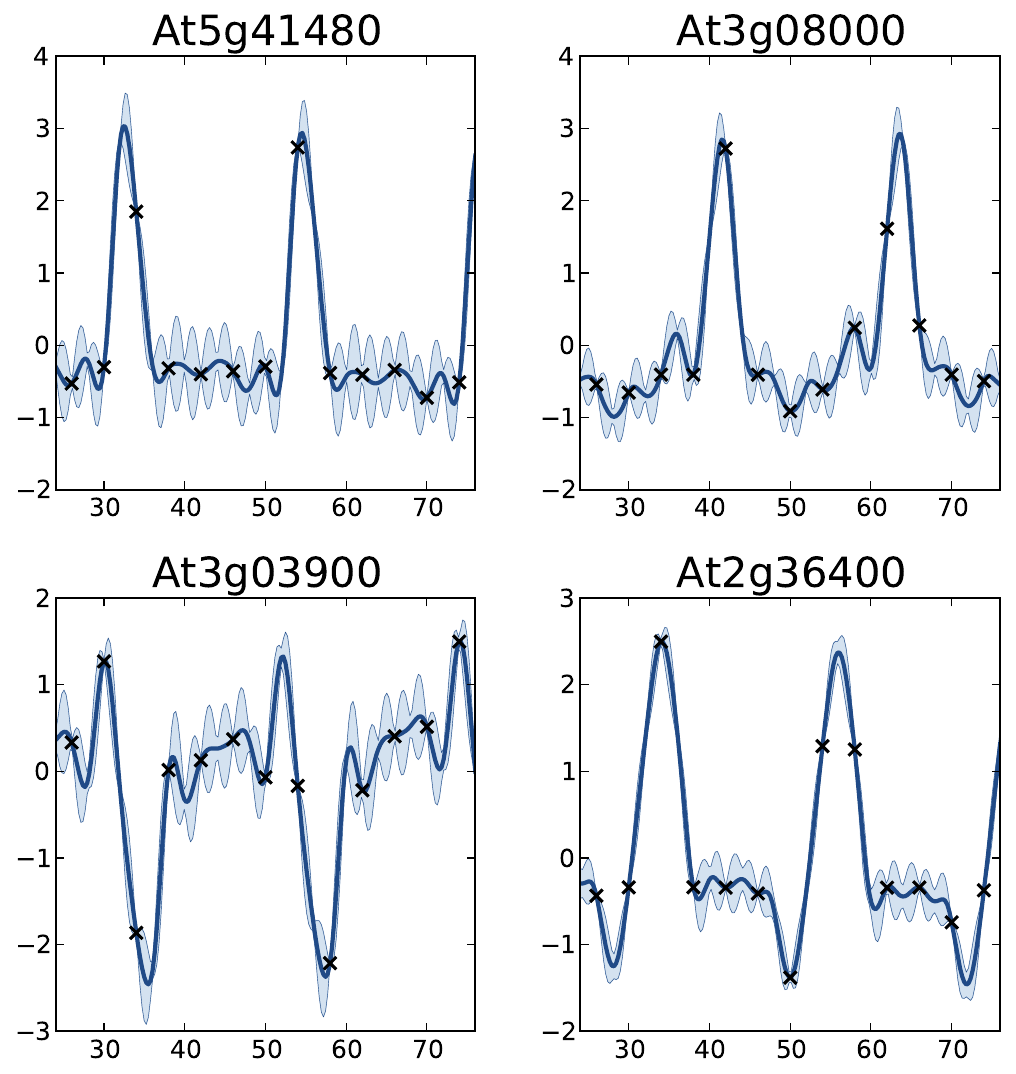}
                \caption{Genes labelled as periodic only by the GP approach.}
        \end{subfigure}
\end{center}
\caption{Examples of genes with different labels. The selected genes
  correspond to the four genes with the highest periodic part
  according to the method that label the gene as periodic. The titles
  of the graphs correspond to the name of the genes (AGI convention).}
\label{fig:compresults}
\end{figure}

\paragraph{}
One interesting comparison between the two methods
is to examine the genes that are classified differently. The available
data from~\cite{edwards2006flowering} allows focusing on the worst classification mistakes made
by one method according to the other. This is illustrated in
Figure~\ref{fig:compresults} which shows the behaviour of the most periodically expressed genes in
$\overline{\mathcal{P}_{GP}}$ according to COSOPT and, conversely, the genes in
$\overline{\mathcal{P}_{COSOPT}}$ with the highest periodicity ratio
$R$. Although it is undeniable that the genes selected only by COSOPT (panel a) present some periodic
component, they also show a strong
non-periodic part, corresponding either to noise or trend. For these genes, the value of the periodicity ratio is: 0.74 (0.10), 0.74
(0.15), 0.63 (0.11), 0.67 (0.05) (means and standard deviations, clockwise from top left) which is close to the classification boundary. On
the other hand, the models suggested only by the GP approach show a
strong periodic signal (we have for all genes $R=1.01\ (0.01)$) with sharp spikes. Another
interesting fact of panel b is that there is at least one observation
associated with each spike which suggests that the
behaviour of the model should not be interpreted as overfitting.

\paragraph{}
This few elements of comparison on a real life case study show some
very promising results, both for the capability of the proposed method
to handle large datasets and for the quality of the
results. Furthermore we believe that the spike shape of the newly
discovered genes may be of particular interest for understanding the
mechanism of the circadian clock. The full results, as well as the
original dataset can be found in the supplementary materials.  

\section{Conclusion}

\paragraph{}
The main purpose of this article is to introduce a new approach for
estimating and extracting the periodic part of a function $f$ given some
observations $f(x_i) = y_i$. As often, the proposed method corresponds to the orthogonal projection
onto a basis of periodic functions. The originality here is to
perform this projection in some RKHS where the partial knowledge given
by the observations can be dealt with elegantly. Previous theoretical results from the mid-1900s allowed us to
derive the expressions of the inner product of RKHS based on Mat\'ern
kernels. Given these results, it was then possible to define a periodic kernel $k_p$ and to decompose $k$ as a sum of sub-kernels $k = k_p + k_a$.

\paragraph{}
We illustrated three fundamental feature of the proposed kernels for GP modelling. First, as we have seen on the benchmark examples, they allow to approximate non-sinusoidal patterns while retaining appropriate filtering of the noise. Second, they provide a natural decomposition of the GP model as a sum of periodic and aperiodic sub-models. Third, they can be reparametrised to
define a wider family of kernel which is of particular interest for decoupling the assumptions on the behaviour of the periodic and aperiodic part of the signal. This approach has proved to increase considerably the prediction ability of the model on the Mauna Loa Observatory dataset.

\paragraph{}
The probabilistic interpretation of the decomposition in sub-models
is of great importance when it comes to define a criterion that quantifies the
periodicity of $f$ while taking into account the uncertainty about it.
This goal was achieved by applying methods commonly used in
GP based sensitivity analysis to define a periodicity ratio. 

\paragraph{}
Although the proposed method can be applied to any time series data, this
work has originally been motivated by the detection of periodically
expressed genes. In practice listing such genes is a key step for a
better understanding of the circadian clock mechanism at a microscopic level. The effectiveness of the method is
illustrated on such data in the last section. The results we obtained are
consistent with the literature but they also feature some new genes
with a strong periodic component. This suggest that the approach
described here is not only theoretically elegant but also efficient in
practice.

\paragraph{}
As a final remark, we would like to stress that the
proposed method is fully compatible with all the features of Gaussian
processes, from the combination of one-dimensional periodic kernels to
obtain periodic kernels in higher dimension
to the use of global optimisation routines such as EGO.

\bigskip
\begin{center}
{\large\bf AKNOWLEDGEMENT}
\end{center}
The authors gratefully acknowledge the support from the BioPreDynProject (Knowledge Based Bio-Economy EU grant Ref 289434) and the BBSRC grant BB/1004769/1.

\bigskip
\begin{center}
{\large\bf SUPPLEMENTARY MATERIAL}
\end{center}

The following datasets are made available under the Public Domain Dedication and License v1.0 whose full text can be found at: \url{http://www.opendatacommons.org/licenses/pddl}.
\begin{description}

\item[Case study dataset:] Original dataset with the gene
  expressions for each gene at each time point. (csv file)

\item[Case study results:] File regrouping the available results
  from~\cite{edwards2006flowering} and the one obtained in the
  application section. For both methods, the file gives the value of
  the criterion and the estimated period. (csv file)

\end{description}

%%%%%%%%%%%%%%%%%%%%%%%%%%%%%%%%%%%%%%%%%%%%%%%%%%%%%%%%%%%%%%%%%%%%%%%%%%%%%%%%%%%%%%%%%%%%%%%%%%%
%%%%%%%%%%%%%%%%%%%%%%%%%%%%%%%%%%%%%%%%%%%%%%%%%%%%%%%%%%%%%%%%%%%%%%%%%%%%%%%%%%%%%%%%%%%%%%%%%%%
%%%%%%%%%%%%%%%%%%%%%%%%%%%%%%%%%%%%%%%%%%%%%%%%%%%%%%%%%%%%%%%%%%%%%%%%%%%%%%%%%%%%%%%%%%%%%%%%%%%
%%%%%%%%%%%%%%%%%%%%%%%%%%%%%%%%%%%%%%%%%%%%%%%%%%%%%%%%%%%%%%%%%%%%%%%%%%%%%%%%%%%%%%%%%%%%%%%%%%%
\appendix
\bigskip
\begin{center}
{\large\bf APPENDIX}
\end{center}
%%%%%%%%%%%%%%%%%%%%%%%%%%%%%%%%%%%%%%%%%%%%%%%%%%%%%%%%%%%%%%%%%%%%%%%%%%%%%%%%%%%%%%%%%%%%%%%%%%%
%%%%%%%%%%%%%%%%%%%%%%%%%%%%%%%%%%%%%%%%%%%%%%%%%%%%%%%%%%%%%%%%%%%%%%%%%%%%%%%%%%%%%%%%%%%%%%%%%%%
\section{Norms in Mat\'ern RKHS}
\label{sec:NORM}

%%%%%%%%%%%%%%%%%%%%%%%%%%%%%%%%%%%%%%%%%%%%%%%%%%%%%%%%%%%%%%%%%%%%%%%%%%%%%%%%%%%%%%%%%%%%%%%%%%%
\subsection{Autoregressive processes and RKHS norms}

\paragraph{}
A process is said to be autoregressive (AR) if the spectral density of the kernel
\begin{equation}
S(\omega) = \frac{1}{2 \pi} \int_\mathds{R} k(t) e^{-i \omega t} \dx \omega
\label{eq:SK}
\end{equation}
can be written as a function of the form
\begin{equation}
S(\omega) = \frac{1}{ \left| \sum_{k=0}^m \alpha_k (i \omega)^k \right|^2}
\label{eq:Sform}
\end{equation}
where the polynomial $\sum_{k=0}^m \alpha_k x^k$ is real with no zeros in the right half of the complex plan~\cite{Doob1953}. Hereafter we assume that $m \geq 1$ and that $\alpha_0,\alpha_m \neq 0$.

\paragraph{}
For such kernels, the inner product of the associated RKHS $\mathcal{H}$ is given by \citet{hajek1962linear,kailath1971rkhs,parzen1961approach}
\begin{equation}
\PS{h,g} = \int_a^b (L_t h)(L_t g) \dx t + 2 \sum_{\substack{0 \leq j,k \leq m-1 \\ j+k \text{ even}}} d_{j,k} h^{(j)}(a) g^{(k)}(a)
\label{eq:thP}
\end{equation}
\begin{equation*}
 \begin{split}
  \text{where } L_t h & = \sum_{k=0}^m \alpha_k h^{(k)}(t) \text{ and }d_{j,k} = \sum_{i=\max(0,j+k+1-n)}^{\min(j,k)} (-1)^{(j-i)} \alpha_i \alpha_{j+k+1-i}.
 \end{split}
\end{equation*}

We show in the next section that the Mat\'ern kernels correspond to autoregressive kernels and, for the usual values of $\nu$, we derive the norm of the associated RKHS.

%%%%%%%%%%%%%%%%%%%%%%%%%%%%%%%%%%%%%%%%%%%%%%%%%%%%%%%%%%%%%%%%%%%%%%%%%%%%%%%%%%%%%%%%%%%%%%%%%%%
\subsection{Application to Mat\'ern kernels}
\paragraph{}
Following the pattern exposed in \citet[p. 542]{Doob1953}, the spectral density of a Mat\'ern kernel (Eq.~\ref{eq:SDmat}) can be written as the density of an AR process when $\nu + 1/2$ is an integer. Indeed, the roots of the polynomial $ \frac{2 \nu}{\theta ^2} + \omega ^2 $ are conjugate pairs so it can be expressed as the squared module of a complex number
\begin{equation}
\begin{split}
  \frac{2 \nu}{\theta ^2} + \omega ^2 & = \Big(\omega + \frac{i \sqrt{2 \nu}}{\theta} \Big) \Big(\omega - \frac{i \sqrt{2 \nu}}{\theta} \Big)  = \Big| \omega + \frac{i \sqrt{2 \nu}}{\theta} \Big| ^2.
\end{split}
\end{equation}
 Multiplying by $i$ and taking the conjugate of the quantity inside the module, we finally obtain a polynomial in $i \omega $ with all roots in the left half of the complex plan:
\begin{equation}
   \frac{2 \nu}{\theta ^2} + \omega^2  = \Big| i \omega + \frac{ \sqrt{2 \nu}}{\theta} \Big|^2
\Rightarrow \left( \frac{2 \nu}{\theta ^2} + \omega^2 \right)^{(\nu+1/2)} = 
\left| \left( \frac{\sqrt{2 \nu}}{\theta} + i \omega \right)^{(\nu+1/2)} \right|^2.
\end{equation}
Plugging this expression into Eq.~\ref{eq:SDmat}, we obtain the desired expression of $S_\nu$:

\begin{equation}
S_\nu(\omega) =  \frac1{ \displaystyle \left| \sqrt{ \frac{ \Gamma (\nu) \theta^{2\nu}}{2 \sigma^2 \sqrt{\pi} \Gamma (\nu + 1/2) (2 \nu)^\nu}} \left( \frac{\sqrt{2 \nu}}{\theta} + i \omega \right)^{(\nu+1/2)} \right|^2 }.
\label{eq:SDmatinvfin}
\end{equation}
Using $\Gamma (\nu) = \frac{(2 \nu - 1)! \sqrt{\pi}}{ 2^{2 \nu -1} (\nu - 1/2)!}$, one can derive the following expression of the coefficients $\alpha_k$:
\begin{equation}
 \alpha_k = \sqrt{ \frac{ (2 \nu - 1)! \nu^\nu}{ \sigma^2 (\nu - 1/2)!^2 2^\nu}} \mathcal{C}^k_{\nu+1/2} \left(\frac{\theta}{\sqrt{2 \nu}} \right)^{k-1/2}.
\end{equation}

\paragraph{}

Theses values of $\alpha_k$ can be plugged into Eq.~\ref{eq:thP} to
obtain the expression of the RKHS inner product. The results for $\nu
\in \{1/2,\ 3/2,\ 5/2 \}$ is given by
Eqs.~\ref{eq:PSk12}-\ref{eq:PSk52}
 in the main body of
the article.

%%%%%%%%%%%%%%%%%%%%%%%%%%%%%%%%%%%%%%%%%%%%%%%%%%%%%%%%%%%%%%%%%%%%%%%%%%%%%%%%%%%%%%%%%%%%%%%%%%%
%%%%%%%%%%%%%%%%%%%%%%%%%%%%%%%%%%%%%%%%%%%%%%%%%%%%%%%%%%%%%%%%%%%%%%%%%%%%%%%%%%%%%%%%%%%%%%%%%%%
%%%%%%%%%%%%%%%%%%%%%%%%%%%%%%%%%%%%%%%%%%%%%%%%%%%%%%%%%%%%%%%%%%%%%%%%%%%%%%%%%%%%%%%%%%%%%%%%%%%
%%%%%%%%%%%%%%%%%%%%%%%%%%%%%%%%%%%%%%%%%%%%%%%%%%%%%%%%%%%%%%%%%%%%%%%%%%%%%%%%%%%%%%%%%%%%%%%%%%%

%\bibliographystyle{natbib}
\bibliographystyle{plainnat}
\bibliography{ref}

\end{document}